\documentclass[a4paper,11pt,draft]{amsart}
\usepackage{amssymb}

\newtheorem{theorem}{Theorem}[section]
\newtheorem{lemma}[theorem]{Lemma}
\newtheorem{prop}[theorem]{Proposition}
\newtheorem{corollary}[theorem]{Corollary}

\theoremstyle{definition}

\newtheorem{example}[theorem]{Example}

\theoremstyle{remark}

\numberwithin{equation}{section}

\newcommand{\DD}{{\mathbb D}}
\newcommand{\OO}{{\mathcal O}}

\newcommand{\NN}{{\mathbb N}}

\newcommand{\RR}{{\mathbb R}}
\newcommand{\CC}{{\mathbb C}}
\newcommand{\ZZ}{{\mathbb Z}}

\newcommand{\eps}{\varepsilon}
\DeclareMathOperator{\re}{Re}
\DeclareMathOperator{\psh}{PSH}
\DeclareMathOperator{\pshc}{PSHC}

\DeclareMathOperator{\dist}{dist}

\renewcommand{\phi}{\varphi}

\begin{document}

%%%%%%%%%%%%%%%%%%%%%%%%%%%%%%%%%%%%%%%%%%%%%%%%%%%%%%%%%%%%%%%%%%%%%%
%%%%%%%%%%%%%%%%%%%%%%%%%%%%%%%%%%%%%%%%%%%%%%%%%%%%%%%%%%%%%%%%%%%%%%
\title[Completeness of the graph]
{The pluripolar hull of the graph of a holomorphic function with polar 
singularities}

\author{Armen Edigarian}

\address{Institute of Mathematics, Jagiellonian University,
Reymonta 4/526, 30-059 Krak\'ow, Poland}
\email{edigaria@im.uj.edu.pl}
\thanks{The first author was supported in part by the KBN grant
No. 5 P03A 033 21. The first author is a fellow of the 
A.~Krzy\.zanowski Foundation (Jagiellonian University)}

\author{Jan Wiegerinck}

\address{Faculty of Mathematics, University of Amsterdam,
Plantage Muidergracht 24, 1018 TV, Amsterdam, The Netherlands}
\email{janwieg@science.uva.nl}
%%%\thanks{Support information for the second author.}

%%%\subjclass{Primary 54C40, 14E20; Secondary 46E25, 20C20}

\date{29 June 2002}

%%%\dedicatory{This paper is dedicated to }

%%%\keywords{Pluripolar hull}

\begin{abstract}
We study the pluripolar hull of the graph of a holomorphic function $f$, 
defined on a domain $D\subset \CC$ outside a polar set $A\subset D$.
This leads to a theorem that describes under what conditions 
$f$ is nowhere extendable over $A$, while the graph of $f$ over 
$D\setminus A$ is {\bf not} complete pluripolar. 
\end{abstract}

\maketitle

%%%%%%%%%%%%%%%%%%%%%%%%%%%%%%%%%%%%%%%%%%%%%%%%%%%%%%%%%%%%%%%%%%%%%%
%%%%%%%%%%%%%%%%%%%%%%%%%%%%%%%%%%%%%%%%%%%%%%%%%%%%%%%%%%%%%%%%%%%%%%
\section{Introduction}

Levenberg, Martin and Poletsky \cite{LMP} have conjectured that
if $f$ is a holomorphic function, which is defined on its maximal
domain of existence $D\subset\CC$, then the graph
\begin{equation}
\Gamma_f=\{\big(z,f(z)\big): z\in D\}
\end{equation}
of $f$ over $D$ is a complete pluripolar subset of $\CC^2$, i.e. there
exists a plurisubharmonic function on $\CC^2$ such that it equals
$-\infty$ precisely on $\Gamma_f$ (see e.g.~\cite{K}).
Supports for this conjecture can be found in 
\cite{LMP, LP, W1, W2, W3, E}.

In the same direction is one of the results of the paper.
\begin{theorem}\label{thm11}
Let $A$ be a closed polar set in a domain $D\subset\CC$ 
and let $f$ be a holomorphic function on $D\setminus A$
(we shall write $f\in\OO(D\setminus A)$). 
Then the set $\Gamma_f\cup (A\times\CC)$ is complete pluripolar in 
$D\times\CC$. Moreover, if $f$ does not extend holomorphically over $A$, 
then the set $\Gamma_f\cup (A'\times\CC)$ is complete pluripolar in
$D\times\CC$, where $A'$ is the set of density points of $A$ in $D$. 
\end{theorem}
%%%%%%%%%%%%%%%%%%%%%%%%%%%%%%%%%%%%%%%%%%%%%%%%%%%%%%%%%%%%%%%%%%%%%%

In case $A'=\varnothing$ (i.e. $A$ is a discrete subset of $D$), 
we see that the graph $\Gamma_f$ of $f$ on $D\setminus A$ is complete
pluripolar in $D\times\CC$. Hence, Theorem~\ref{thm11} is a generalization
of the main result of \cite{W2}.

Nevertheless, in \cite{EW} we have shown that the conjecture is false. 
Moreover, for any non-empty set $A'$ there exists a holomorphic 
function defined on $D\setminus A$,
which is nowhere extendable over $A$,
such that the graph $\Gamma_f$ of $f$ on $D\setminus A$ is not
complete pluripolar in $D\times\CC$. The main purpose
of the paper is to clarify this phenomenon.

To study completeness  of the graph we need more definitions. 
Following \cite{Z} (see also \cite{LP}) we introduce for a pluripolar set $E$ 
in a domain $D$ its {\it pluripolar hull in $D$ \/}  as
\begin{equation}
E_D^\ast:=\{z\in D: \text{ for all } u\in\psh(D): u|_E=-\infty
\implies u(z)=-\infty\}.
\end{equation}
It is well-known that a polar set $P$ in $\CC$ is complete polar if
and only if $P$ is of $G_\delta$-type \cite{L}.
In higher dimension, the situation is more complicated.
It is easy to see that any complete pluripolar set is of $G_\delta$-type.
The following result of A.~Zeriahi is a useful counterpart of 
the one-dimensional case.
%%%%%%%%%%%%%%%%%%%%%%%%%%%%%%%%%%%%%%%%%%%%%%%%%%%%%%%%%%%%%%%%%%%%%%
\begin{theorem}[see \cite{Z}]\label{Zer-thm}
Let $E$ be a pluripolar subset of a pseudoconvex domain $D\subset\CC^n$.
If $E_D^\ast=E$ and $E$ is of $G_\delta$-  and $F_\sigma$-type, then
$E$ is complete pluripolar in $D$.
\end{theorem}
%%%%%%%%%%%%%%%%%%%%%%%%%%%%%%%%%%%%%%%%%%%%%%%%%%%%%%%%%%%%%%%%%%%%%%

The next result is the main one of this paper. 
It is a one-variable necessary and sufficient condition on holomorphic 
functions on a domain in $\CC$ that guarantees that the graph of such a 
function is complete pluripolar.
%%%%%%%%%%%%%%%%%%%%%%%%%%%%%%%%%%%%%%%%%%%%%%%%%%%%%%%%%%%%%%%%%%%%%%
%
%%%%%%%%%%%%%%%%%%%%%%%%%%%%%%%%%%%%%%%%%%%%%%%%%%%%%%%%%%%%%%%%%%%%%%
\begin{theorem}\label{main-thm}
Let $D$ be a domain in $\CC$ and let $A$ be a closed polar subset of $D$.
Suppose that $f\in\OO(D\setminus A)$ and that $z_0\in A$.
Then the following conditions are equivalent:
\begin{enumerate}
\item 
$(\{z_0\}\times\CC)\cap (\Gamma_f)^\ast_{D\times\CC}=\varnothing$;
\item the set $\{z\in D\setminus A: |f(z)|\ge R\}$
is not thin at $z_0$ for any $R>0$.
\end{enumerate}
Moreover, if the set $\{z\in D\setminus A: |f(z)|\ge R\}$
is thin at $z_0$ for some $R>0$, then there exists a $w_0\in\CC$, $|w_0|\le R$, such 
that $\Gamma_f$ is not pluri-thin at $(z_0,w_0)$ (and therefore 
$(z_0,w_0)\in(\Gamma_f)^\ast_{D\times\CC}$).
\end{theorem}
%%%%%%%%%%%%%%%%%%%%%%%%%%%%%%%%%%%%%%%%%%%%%%%%%%%%%%%%%%%%%%%%%%%%%%
%
%%%%%%%%%%%%%%%%%%%%%%%%%%%%%%%%%%%%%%%%%%%%%%%%%%%%%%%%%%%%%%%%%%%%%%
For the definition of thinness and pluri-thinness see \eqref{thin} below and 
cp.~\cite{Brel,K,R}.

Since any graph is of $G_\delta$-  and $F_\sigma$-type, as a corollary
of Theorems~\ref{Zer-thm},~\ref{main-thm} we obtain the following.

\begin{corollary}
Let $D$ be a domain in $\CC$ and let $A$ be a closed polar subset of $D$.
Suppose that $f\in\OO(D\setminus A)$. Then $\Gamma_f$
is complete pluripolar in $D\times \CC$ if and only if 
for all $R>0$ the set
$\{z\in D\setminus A: |f(z)|\ge R\}$ is not thin at any point of A.
\end{corollary}

By Theorem~\ref{thm11} (see also \cite{W2}) in case $A$ has no limit point in 
$D$, $\Gamma_f$ is complete pluripolar (and in particular the  regularity 
condition is satisfied) for every $f\in\OO(D\setminus A)$, which does not
extend holomorphically over $A$.

For example, the corollary applies to the function $f(z)=\exp (1/z)$, which
is in $\OO(\CC\setminus\{0\})$ and does not extend holomorphically to $0$. 
Moreover, $\{z\in\CC\setminus\{0\}: \re(1/z)\ge \ln R\}$ is not thin at $0$
for any $R>0$. Hence, $\Gamma_f$ is complete pluripolar in $\CC^2$ as was
already proven in \cite{W1}, \cite{W2}.
More interestingly, consider the function $f(z)=1/\sin(\pi/z)$ in
$\OO(\CC\setminus{A})$, where $A=\{0\}\cup\{1/n:\ n\in\ZZ\}$. It is an amusing
exercise to verify that in this case for any $R>0$ the set
$\{z\in\CC\setminus A: |f(z)|\ge R\}$ is not thin at any point of $A$
showing again that $\Gamma_f$ is complete pluripolar in  $\CC^2$.

%Recall that in case of isolated singularity the condition (1) in
%Theorem~\ref{main-thm} is not fullfiled only for removable singularity
%(see \cite{W2}). Hence, we obtain the following generalization
%of the well-known result.
%
%\begin{corollary} Let $D$ be a domain in $\CC$ and let 
%$f\in\OO(D\setminus\{z_0\})$, where $z_0\in D$. Then
%$f$ has removable singularity if and only if there exists an $R>0$
%such that the set $\{z\in D\setminus\{z_0\}: |f(z)|\ge R\}$ is
%thin at $z_0$.
%\end{corollary}

The first named author thanks Marek Jarnicki and Peter Pflug for
very helpful remarks. A conversation we had with W\l odzimierz Zwonek
was very useful too.
%%%%%%%%%%%%%%%%%%%%%%%%%%%%%%%%%%%%%%%%%%%%%%%%%%%%%%%%%%%%%%%%%%%%%%

\section{Pluripolar completeness}

Let $\Omega$ be a domain in $\CC^n$. We say that a plurisubharmonic function
$u:\Omega\to[-\infty,\infty)$ is continuous if $e^u:\Omega\to [0,\infty)$ 
is a continuous function and write $u\in\pshc(\Omega)$.

The hart of this section is Theorem \ref{thm1}.
Together with Corollary \ref{cor-11} and Theorem~\ref{thm9}
it leads easily  to Corollary \ref{thm10}.
This in turn forms a basic ingredient in the proof of Theorem~\ref{main-thm}.

\begin{theorem}[cf.~\cite{W2}]\label{thm1}
Let $K$ be a compact polar set in a domain $D\subset\CC$ 
and let $f\in\OO(D\setminus K)$. 
Then there exists a continuous plurisubharmonic function 
$u\in\pshc(D\times\CC)$ such that
\begin{equation*}
\{(z,w)\in D\times\CC: u(z,w)=-\infty\}=
\{(z,f(z)):z\in D\setminus K\}\bigcup (K\times\CC).
\end{equation*}
\end{theorem}

We postpone the proof Theorem~\ref{thm1} and start with some auxiliary results 
leading to Proposition \ref{lem3}.

Put  $\DD_r(z_0):=\{z\in\CC: |z-z_0|<r\}$, $r>0$,
$\DD_r:=\DD_r(0)$, and $\DD:=\DD_1$.

\begin{lemma}\label{lem1} Let $K$ be a compact set in the unit disk $\DD$
and let $f$ be a holomorphic function in $\DD\setminus K$. Then there exist
$f_1\in\OO(\DD)$ and $f_2\in\OO(\CC\setminus K)$ such that
\begin{equation}
f=f_1+f_2\qquad\text{ on } \DD\setminus K
\end{equation}
and $\lim_{|z|\to+\infty}f_2(z)=0$.
\end{lemma}

\begin{proof}%[Proof of Lemma \ref{lem1}]
Take $r\in(0,1)$ such that 
$f$ is holomorphic in the annulus
\begin{equation*}
A(r):=\{z\in\CC:r<|z|<1\}.
\end{equation*}
Then 
\begin{equation}
f(z)=\sum_{k=-\infty}^\infty a_k z^k,\quad z\in A(r).
\end{equation}
Take $\tilde f_1(z):=\sum_{k=0}^\infty a_k z^k$ and
$\tilde f_2(z):=\sum_{k=1}^\infty a_{-k} z^{-k}$.
Then $\tilde f_1$ extends holomorphically to $\DD$ 
(we denote the extension by $f_1$)
and
$\tilde f_2$ extends holomorphically to $\CC\setminus \overline{\DD}_r$
(we denote the extension by $f_2$).
We extend $f_2$ to $\overline{\DD}_r\setminus K$ by the formula
\begin{equation*}
f_2(z)=f(z)-f_1(z),\quad z\in \overline{\DD}_r\setminus K.
\end{equation*}
Note that $f_2\in\OO(\CC\setminus K)$.
\end{proof}

\begin{lemma}\label{lem2} 
Let $A$ be a polar set in $\DD$.
Then there exist a sequence $r_n\nearrow1$ such that 
\begin{equation*}
A\cap \{z\in\DD:|z|=r_n\}=\varnothing, \quad n\ge1.
\end{equation*}
\end{lemma}

\begin{proof}%[Proof of Lemma \ref{lem2}]
This easily follows from the fact that linear measure
(i.e., 1-di\-men\-sional Hausdorff measure) of $A$ is zero
(see \cite{T}, Theorem III.19).
\end{proof}

\begin{prop}\label{lem3}
Let $K$ be a compact polar set in a domain $D\subset\CC$
and let $f\in\OO(D\setminus K)$. Then there exist
$f_1\in\OO(D)$ and $f_2\in\OO(\CC\setminus K)$ such that
\begin{equation}
f=f_1+f_2\qquad\text{ on } D\setminus K
\end{equation}
and $\lim_{|z|\to\infty}f_2(z)=0$.
\end{prop}

\begin{proof}%[Proof of Proposition~\ref{lem3}]
For any $z_0\in K$ there exists $r_0>0$ such that
$\DD_{r_0}(z_0)\subset D$ and $\partial\DD_{r_0}(z_0)\cap K=\varnothing$
(use Lemma~\ref{lem2}). 
From the compactness of $K$ it follows that
there exist points $z_1,\dots,z_N\in K$ and $r_1,\dots,r_N>0$ such that

a) $K\subset\bigcup_{j=1}^N \DD_{r_j}(z_j)$,

b) $\partial\DD_{r_j}(z_j)\cap K=\varnothing$, $j=1,\dots,N$.

Now, we construct functions $g_1,\dots,g_N\in\OO(\CC\setminus K)$
such that $f-(g_1+\dots+g_N)\in\OO(D)$.

By Lemma~\ref{lem1} there exists $g_1\in\OO(\CC\setminus K_1)$,
where $K_1=K\cap \DD_{r_1}(z_1)$, such that $f-g_1\in\OO(\DD_{r_1}(z_1))$.

Again by Lemma~\ref{lem1} there exists $g_2\in\OO(\CC\setminus K_2)$,
where $K_2=(K\cap \DD_{r_2}(z_2))\setminus K_1$, such that
$f-g_1-g_2\in\OO(\DD_{r_1}(z_1)\bigcup \DD_{r_2}(z_2))$.
Repeating this process, we construct $g_1,\dots,g_N\in\OO(\CC\setminus K)$
such that $f-(g_1+\dots+g_N)\in\OO(\bigcup_{j=1}^N \DD_{r_j}(z_j))$.
\end{proof}

For the proof of Theorem~\ref{thm1} we need the following result.
The use of Jacobi series in this setting is well-known, 
cf.~\cite{S1}, \cite{S2}.

\begin{prop}%[cf.~\cite{S1}, \cite{S2}]
\label{app-thm} 
Let $K$ be a compact polar set in $\CC$ and let 
$f\in\OO(\CC\setminus K)$. Then there exist rational functions
$r_n(z)=\frac{p_n(z)}{q_n(z)}$ of degree $\le n$
(i.e., $\deg p_n\le n$ and $\deg q_n\le n$)
 with singularities in $K$ 
(i.e., the zeros of $q_n$ lie in $K$)
such that for any compact set $L\subset\CC\setminus K$
we have
\begin{equation}
\|f-r_n\|_L^{\frac 1n}\to0,\quad \text{if $n\to\infty$}.
\end{equation}
Here $\|\cdot\|_L$ means $\sup$-norm in $L$.
\end{prop}

\begin{proof}[Proof of Proposition~\ref{app-thm}]
According to Proposition~\ref{lem3} there exist $f_1\in\OO(\CC)$ and 
$f_2\in\OO(\CC\setminus K)$ such that $f=f_1+f_2$ on $\CC\setminus K$ and
$f_2(z)\to0$ when $|z|\to\infty$. So, we may assume that
$f(z)\to0$ when $|z|\to\infty$.

Let $U_n=\{z\in\CC:\dist(z,K)<\frac1{n}\}$, $n=1,2,\ldots$. Then 
$U_{n+1}\Subset U_n$ and $\cap_{n=1}^\infty U_n=K$.

Let $q_m$ denote the $m$-th Fekete polynomial (see e.g. \cite{R}). Since $K$ is polar,
$\|q_m\|_K^{\frac1m}\to0$, $m\to\infty$.

Take $n\in\NN$ so large that $L\Subset\CC\setminus U_n$ and next
 $m(n)\in\NN$ such that 
$\|q_m\|_K^{\frac1m}<\frac1{n^4}$ for  $m\ge m(n)$. Then for $z\in\CC\setminus
U_{n^2}$
\begin{equation*}
|q_m(z)|=\prod_{i=1}^m|z-z_i|\ge n^{-2m}> 
n^{2m}\|q_m\|_K=:\rho_m,
\end{equation*}
while for $z\in\CC\setminus
U_{n}$
\begin{equation*}
|q_m(z)|\ge n^{-m}>n^m\rho_m.
\end{equation*}

Put $V_m=\{z\in\CC: |q_m(z)|> \rho_m\}$. Then
$L\subset U_n^c\subset U_{n^2}^c\subset V_m$. 
Set
\begin{multline*}
F_m(z,w)=\frac{1}{2\pi i q_m(z)}\int_{\partial V_m}
\frac{f(\zeta)}{wq_m(\zeta)-1}\cdot\frac{q_m(\zeta)-q_m(z)}{\zeta-z}d\zeta,\\
z\in V_m, |w|<\frac1{\rho_m}.
\end{multline*}
Then $F_m\in\OO\big(V_m\times\DD(\frac1{\rho_m})\big)$ and 
\[
F_m\Big(z,\frac{1}{q_m(z)}\Big)=
\frac1{2\pi i}\int_{\partial V_m}\frac{f(\zeta)}{\zeta-z}
d\zeta=f(z),
\]
for $z\in V_m$, because $f$ vanishes at $\infty$. Put
\begin{equation*}
c_{mk}(z)=-\frac{1}{2\pi i}\int_{\partial V_m}
f(\zeta)q_m^k(\zeta)\frac{q_m(\zeta)-q_m(z)}
{\zeta-z}d\zeta.
\end{equation*}
Observe that $c_{mk}$ is a polynomial of degree $\le m-1$ and that
$$F(z,w)=\sum_{k=0}^\infty \frac{c_{mk}(z)}{q_m(z)}w^k.$$ Therefore,
\begin{equation*}
f(z)=\sum_{k=0}^\infty\frac{c_{mk}(z)}{q_m^{k+1}(z)}.
\end{equation*}
Put $f_N(z)=\sum_{k=0}^{N-1}\frac{c_{mk}(z)}{q_m^{k+1}(z)}$. This is a rational
function of degree $mN$. Let us
estimate $f-f_N$ on  $L\Subset U_n^c$.
%$\{z\in\CC:|q_m(z)|>\big(\frac{1}{n}\big)^m\|q_m\|_K\}$.
We have
\begin{equation*}
|f(z)-f_N(z)|\le \sum_{k=N}^\infty M_{1}\frac{\rho_{m}^k}
{\min_L |q_m(z)|^{k+1}}\le
M_{2}\sum_{k=N}^\infty\big(\frac{n^{2m}}{n^{3m}}\big)^{k}\le
M_{3}\big(\frac{1}{n}\big)^{mN}.
\end{equation*}
where $M_{j}>0$ may depend on $n,m$ but
are independent of $N$. Hence, for sufficiently big $N$ we have
$|f(z)-f_N(z)|^{\frac1{mN}}\le \frac{1}{n}$.
\end{proof}

\begin{proof}[Proof of Theorem~\ref{thm1}]
In view of Proposition~\ref{lem3} we may assume that $D=\CC$.

By Theorem~\ref{app-thm} there exist rational functions
$r_n=\frac{p_n}{q_n}$ of degree $n$, $n\ge1$,
such that
\begin{equation*}
\|f-r_n\|_L^{\frac 1n}\to0
\end{equation*}
for any compact set $L\subset\CC\setminus K$.
We may assume that 
$q_n(z)=(z-a_{n1})\dots (z-a_{nm_n})$, where 
$a_{n1},\dots,a_{m_nn}\in K$, $m_n\le n$.

Put
\begin{equation*}
h_n(z,w):=\frac1n\log\left|\Big(w-r_n(z)\Big)q_n(z)\right|.
\end{equation*}
Then $h_n$ is a continuous plurisubharmonic function in $\CC\times\CC$.

\medskip

We may assume that $K\subset\DD$. For any $\nu\ge2$ we put
\begin{equation*}
D_\nu=\{z\in\DD_\nu:\dist(z,K)>\frac{1}{\nu}\}.
\end{equation*}
Note that $D_\nu\subset\subset D_{\nu+1}$
and $\bigcup_{\nu} D_{\nu}=\CC\setminus K$. 

For $z\in D_{\nu}$ and $n\in\NN$ we have
\begin{multline}
h_n(z,f(z))=\frac{1}{n}\log|(f(z)-r_n(z))q_n(z)|=
\frac{1}{n}\log|f(z)-r_n(z)|\\
+\frac{1}{n}\log|q_n(z)|\le\log\|f-r_n\|_{D_\nu}^{\frac1n}+\log (\nu+1)
\\=\log\left( (\nu+1)\|f-r_n\|_{D_\nu}^{\frac1n}\right).
\end{multline}
Recall that $\|f-r_n\|_{D_\nu}^{\frac1n}\to0$, when $n\to0$. Hence, for
any $\nu\in\NN$ there exists $n_1(\nu)$ such that
for any $n\ge n_1(\nu)$
\begin{equation}\label{eq2.6}
h_n(z,f(z))\le -\nu,\quad z\in D_\nu.
\end{equation}

\medskip

Fix $\nu\in\NN$. For $z\in D_\nu$, $w\in\DD_\nu$,
with $|w-f(z)|>\frac{1}{\nu}$ and $n\ge n_1(\nu)$ we have, 
using \eqref{eq2.6}
\begin{multline}
h_n(z,w)=\frac{1}{n}\log|(w-f(z))q_n(z)+(f(z)-r_n(z))q_n(z)|\ge\\
\frac{1}{n}\log\Big|\left(\frac{1}{\nu}\right)^{n+1}-e^{-\nu n}\Big|\ge
\frac{1}{n}\log\frac12\Big(\frac1\nu)^{n+1}.
\end{multline}
Hence, there exists $n_2(\nu)\ge n_1(\nu)$ such that for
$n\ge n_2(\nu)$ we have
\begin{equation}\label{eq2.8}
h_n(z,w)\ge-\log(\nu+1).
\end{equation}

\medskip

Let us estimate $h_n$ on $\DD_\nu\times\DD_\nu$.
Since $\|f-r_n\|_{\{z\in\CC:|z|=\nu\}}^{\frac1n}\to0$, 
for sufficiently big $n\in\NN$ we have
\begin{equation*}
\max_{|z|\le\nu}|p_n(z)|=\max_{|z|=\nu}|p_n(z)|\le 
C_\nu\max_{|z|=\nu}|q_n(z)|,
\end{equation*}
where $C_\nu=\max_{|z|=\nu}|f(z)|+1$.
Hence, there exists $n_3(\nu)\ge n_2(\nu)$ such that
\begin{multline}
h_n(z,w)=\frac1n\log|wq_n(z)-p_n(z)|\\
\le\frac1n\log\Big(\nu(\nu+1)^n+C_\nu(\nu+1)^n\Big)\le \log(\nu+2),\\
n\ge n_3(\nu), (z,w)\in\DD_\nu\times\DD_\nu.
\end{multline}

For any $\nu$ fix an $n(\nu)\ge n_3(\nu)$.
Now, we consider the plurisubharmonic functions
\begin{equation*}
u_\nu(z,w)=\max\{h_{\nu(n)}(z,w)-\log(\nu+2), -\nu-\log(\nu+2)\}.
\end{equation*}
As $u_\nu$ is negative on $\DD_\nu\times\DD_\nu$, the series
\begin{equation}
u(z,w):=\sum_{\nu=2}^\infty \frac{1}{\nu^2} u_\nu(z,w).
\end{equation}
represents a plurisubharmonic function. It is not identically $-\infty$
because of \eqref{eq2.8}, while \eqref{eq2.6} shows that it is $-\infty$ 
on the graph of $f$.
Finally, the function $u(z,w)+H(z)$, where $H(z)$ is the
Evans' potential of $K$ (see e.g.~\cite{R}), satisfies all our conditions.
\end{proof}

%%%%%%%%%%%%%%%%%%%%%%%%%%%%%%%%%%%%%%%%%%%%%%%%%%%%%%%%%%%%%%%%%%%%%%
\section{Results from potential and pluripotential theories}
%%%%%%%%%%%%%%%%%%%%%%%%%%%%%%%%%%%%%%%%%%%%%%%%%%%%%%%%%%%%%%%%%%%%%%
In  \ref{ss31} we assemble some results mainly of Levenberg and Poletsky,
\cite{LP}.
Next, in \ref{ss32} we turn to potential theory in $\CC$. We will obtain some
auxiliary results that look very classical. Yet we could not find them in the 
literature.

\subsection{Pluripotential theory} \label{ss31}

 For a pluripolar set $E$ in a domain $D$ in $\CC^n$,
 Levenberg and Poletsky define
the {\it negative pluripolar hull of $E$\/} in $D$ as
\begin{equation*}
E_D^{-}:=\{z\in D: \text{ for all } u\in\psh(D), u\le0:u|_E=-\infty
\implies u(z)=-\infty\}.
\end{equation*}

We have the following relation between pluripolar and
negative pluripolar hulls.
%%%%%%%%%%%%%%%%%%%%%%%%%%%%%%%%%%%%%%%%%%%%%%%%%%%%%%%%%%%%%%%%%%%%%%
\begin{theorem}[see \cite{LP}]\label{lev-pol-thm}
Let $D$ be a pseudoconvex domain in $\CC^n$
and let $E\subset D$ be pluripolar. Suppose that
$D=\bigcup_{j} D_j$, where $D_j\subset D_{j+1}$ form an increasing
sequence of relatively compact pseudoconvex subdomains of $D$. Then
\begin{equation*}
E_D^\ast=\bigcup_j (E\cap D_j)_{D_j}^{-}.
\end{equation*}
\end{theorem}
%%%%%%%%%%%%%%%%%%%%%%%%%%%%%%%%%%%%%%%%%%%%%%%%%%%%%%%%%%%%%%%%%%%%%%
As an immeadiate corollary of Theorem~\ref{lev-pol-thm} we get the following.
\begin{corollary}\label{cor-11}
Let $G$ be a domain in $\CC$ and let $E\subset G\times\CC$ be a pluripolar set,
which is of $F_\sigma$- and $G_\delta$-type. 
Then $E$ is complete pluripolar in  $G\times\CC$ iff $E\cap (D\times\CC)$ is
complete pluripolar in $D\times\CC$ for any domain $D\Subset G$.
\end{corollary}
%%%%%%%%%%%%%%%%%%%%%%%%%%%%%%%%%%%%%%%%%%%%%%%%%%%%%%%%%%%%%%%%%%%%%%
\begin{proof} Note that if $E$ is complete pluripolar in $G\times\CC$
then $E$ is complete pluripolar in any smaller domain.

Take a sequence of domains $D_1\Subset D_2\Subset\dots\Subset G$ such that
$G=\cup_{\nu=1}^\infty D_\nu$. By Theorem~\ref{Zer-thm} 
$E$ is complete pluripolar in $G\times\CC$ iff $E^\ast_{G\times\CC}=E$.
By Theorem~\ref{lev-pol-thm} we have
\begin{equation*}
E^\ast_{G\times\CC}=\cup_{\nu=1}^\infty
\big(E\cap (D_\nu\times\DD_\nu)\big)^{-}_{D_\nu\times\DD_\nu}.
\end{equation*}
We know that
\begin{equation*}
\big(E\cap(D_{\nu+1}\times\CC)\big)^\ast_{D_{\nu+1}\times\CC}\supset
\big(E\cap(D_{\nu}\times\DD_\nu)\big)^{-}_{D_{\nu}\times\DD_\nu},
\quad \nu=1,2,\dots
\end{equation*}
So, if $\big(E\cap(D_{\nu+1}\times\CC)\big)^\ast_{D_{\nu+1}\times\CC}=
E\cap(D_{\nu+1}\times\CC)$, then $E^\ast_{G\times\CC}=E$.
\end{proof}

%%%%%%%%%%%%%%%%%%%%%%%%%%%%%%%%%%%%%%%%%%%%%%%%%%%%%%%%%%%%%%%%%%%%%%
In view of Theorem~\ref{thm1} recall the following result of
A.~Zeriahi (see \cite{Z}).

\begin{theorem}\label{thm9}
Let $D$ be a pseudoconvex domain in $\CC^n$ and let 
$Z$ be a closed complete pluripolar set in $D$. Then
there exists a continuous plurisubharmonic function 
$u\in\pshc(D)$ such that 
$\{z\in D: u(z)=-\infty\}=Z$.
\end{theorem}

As a corollary of Theorems~\ref{thm1} and~\ref{thm9}, Corollary~\ref{cor-11} 
we get the following.
\begin{corollary}\label{thm10}
Let $A$ be a closed polar set in a domain $D\subset\CC$ 
and let $f\in\OO(D\setminus A)$. 
Then there exists a continuous plurisubharmonic function
$u\in\pshc(D\times\CC)$ such that
\begin{equation*}
\{(z,w)\in D\times\CC: u(z,w)=-\infty\}=
\{(z,f(z)):z\in D\setminus A\}\bigcup (A\times\CC).
\end{equation*}
In particular,
$(\Gamma_f)^\ast_{D\times\CC}\cap \big[(D\setminus A)\times \CC\big]=
\Gamma_f\cap \big[(D\setminus A)\times\CC\big]$.
\end{corollary}
%%%%%%%%%%%%%%%%%%%%%%%%%%%%%%%%%%%%%%%%%%%%%%%%%%%%%%%%%%%%%%%%%%%%%%
\begin{proof} Let $G\Subset D$ be a domain such that 
$\partial G\cap A=\varnothing$. Then $A\cap G$ is a compact
polar subset of $G$. By Theorem~\ref{thm1} the set
$\{(z,f(z)): z\in G\setminus A\}\cup \big((A\cap G)\times\CC\big)$ is
complete pluripolar in $G\times\CC$. Using Corollary~\ref{cor-11}
we get that $\{(z,f(z)):z\in D\setminus A\}\bigcup (A\times\CC)$
is complete pluripolar in $D\times\CC$. The existence of a 
continuous plurisubharmonic function follows from Theorem~\ref{thm9}.
\end{proof}
%%%%%%%%%%%%%%%%%%%%%%%%%%%%%%%%%%%%%%%%%%%%%%%%%%%%%%%%%%%%%%%%%%%%%%
Let $E$ be a subset of a domain $D\subset\CC^n$. The
{\it pluriharmonic measure\/} of $E$ relative to $D$
(or, {\it relative extremal function\/}) is defined as follows
\begin{multline}\label{eq-ext}
\omega(z,E,D)=-\sup\{u(z):u\in\psh(D)\text{ and }u\le-1\text{ on }E,
u\le0\text{ on }D\},\\
z\in D.
\end{multline}
We refer to \cite{K} for details. Note that
relative extremal function can be defined as in \eqref{eq-ext} for 
an open set $D$. In that case we have 
$\omega(z,E,D)=\omega(z,E\cap D_j, D_j)$, $z\in D_j$, where
$D=\cup_{j} D_j$, $D_j$ are disjoint and connected.

Recall the following important result.
\begin{theorem}[see \cite{LP}]\label{lev-pol2-thm}
Let $D$ be a domain in $\CC^n$ and let $E$ be
a pluripolar set in $D$. Then
\begin{equation*}
E_D^-=\{z\in D: \omega(z,E,D)>0\}.
\end{equation*}
\end{theorem}
%%%%%%%%%%%%%%%%%%%%%%%%%%%%%%%%%%%%%%%%%%%%%%%%%%%%%%%%%%%%%%%%%%%%%%
%
%%%%%%%%%%%%%%%%%%%%%%%%%%%%%%%%%%%%%%%%%%%%%%%%%%%%%%%%%%%%%%%%%%%%%%
The next localization principle is an important tool in our theory. 
A very special case of this principle can be found in \cite{W2}.
%%%%%%%%%%%%%%%%%%%%%%%%%%%%%%%%%%%%%%%%%%%%%%%%%%%%%%%%%%%%%%%%%%%%%%
\begin{theorem}\label{thm-loc} Let $G$ be a domain in $\CC^n$ and let $E$
be a pluripolar subset of $G$. Assume that $u$ is a
plurisubharmonic function on $G$ such that
\begin{equation}
Z=\{z\in G: u(z)=-\infty\}\supset E
\end{equation}
and $u$ is locally bounded on $G\setminus Z$ (e.g. $u\in\pshc(G)$).
Then for any domain $D\Subset G$ and any neighborhood
$U\supset Z$ we have
\begin{equation}
\omega(z,E\cap D,D)=\omega(z,E\cap D,U\cap D),
\quad z\in U\cap D.
\end{equation}
\end{theorem}
%%%%%%%%%%%%%%%%%%%%%%%%%%%%%%%%%%%%%%%%%%%%%%%%%%%%%%%%%%%%%%%%%%%%%%

\begin{proof} Fix a domain $ D\Subset G$.
Since $U\cap  D\subset D$, we have the inequality ''$\ge$''.

Let us show the inequality ''$\le$''. Note that $Z\cap\overline{ D}$
is a compact set. Assume that $h\in\psh( D\cap U)$ is such that
$h\le-1$ on $E\cap D$ and $h\le0$ on $U\cap D$.
Fix $\epsilon>0$. Then there exist $\alpha>0$, $\beta\in\RR$ such that
$\alpha u+\beta<0$ on $ D$, $\alpha u+\beta\ge-\epsilon$ on
$(\partial U)\cap\overline{ D}$. 

Consider the function
\begin{equation*}
v_\epsilon(z):=\begin{cases}
\max\{h(z)-\epsilon,\alpha u(z)+\beta\} &\text{if }z\in D\cap U,\\
\alpha u(z)+\beta &\text{if } z\in D\setminus U.
\end{cases}
\end{equation*}
Note that $v_\epsilon\in\psh( D)$ is a negative plurisubharmonic function
such that $v_\epsilon\le-1$ on $E\cap D$. Hence,
\begin{equation*}
-\omega(z,E\cap D, D)\ge v_\epsilon(z)\ge h(z)-\epsilon,
\quad z\in  D\cap U.
\end{equation*}
Let $\epsilon\to0$. Then
$-\omega(z,E\cap D, D)\ge h(z)$, $z\in D\cap U$. And, therefore,
\begin{equation*}
\omega(z,E\cap D, D)\le \omega(z,E\cap D,U\cap D).
\end{equation*}
\end{proof}

%%%%%%%%%%%%%%%%%%%%%%%%%%%%%%%%%%%%%%%%%%%%%%%%%%%%%%%%%%%%%%%%%%%%%%
\subsection{Potential theory}\label{ss32}
%%%%%%%%%%%%%%%%%%%%%%%%%%%%%%%%%%%%%%%%%%%%%%%%%%%%%%%%%%%%%%%%%%%%%%
%
%%%%%%%%%%%%%%%%%%%%%%%%%%%%%%%%%%%%%%%%%%%%%%%%%%%%%%%%%%%%%%%%%%%%%%
The following result from potential theory is very important in our
considerations. Since we could not find the proof in the literature, 
we give the proof.
\begin{theorem}\label{thm-imp7}
Let $D$ be a bounded domain in $\CC$ and let $S\subset D$ be a closed disc.
Assume that $K\subset\partial D$ is a compact polar set. Then
\begin{equation*}
\omega(z,S,D)=\inf\{\omega(z,S,D\cup U): K\subset U\text{ open}\},\quad z\in D.
\end{equation*}
Moreover, if $z\in K$ is a regular  boundary point of $D$ 
then
\begin{equation*}
\inf\{\omega(z,S,D\cup U): K\subset U\text{ open}\}=0.
\end{equation*}
\end{theorem}
%%%%%%%%%%%%%%%%%%%%%%%%%%%%%%%%%%%%%%%%%%%%%%%%%%%%%%%%%%%%%%%%%%%%%%
\begin{proof}
To prove the first statement, it is clear that for each $z\in D$
\begin{equation*}
\omega(z,S,D)\le\inf\{\omega(z,S,D\cup U): K\subset U\text{ open}\}.
\end{equation*}
Next, because $K$ is a compact polar set, there exists
a subharmonic function $h\le 0$ that is defined in a neighborhood of
$\bar D$ and such that  $h(w)=-\infty$ if and only if $w\in K$.
Let $z\in D$ and $\eps>0$. Let $U_0$ be a small neighborhood of $K$, such that
$\eps h<-1$ on $U_0$. Then
\begin{multline*}
\inf\{\omega(z,S,D\cup U): K\subset U\text{ open}\} +\eps h(z) \le \\
\omega(z,S,D\cup U_0) +\eps h(z)\le \omega(z,S,D),
\end{multline*}
because $\omega(z,S,D\cup U_0) +\eps h(z)$ is a subharmonic function on
$D\setminus S$ with boundary values less than the boundary values of
$\omega(z,S,D)$.
Since $h(z)$ is a fixed finite value, and $\eps$ is arbitrary $>0$, we have
\begin{equation*}
\omega(z,S,D)\ge\inf\{\omega(z,S,D\cup U): K\subset U\text{ open}\}.
\end{equation*}
and the first statement is established.

For the second statement we first observe that for every $\eps>0$ the set
\begin{equation*}
F^\eps_U=\{ z\in \overline{D\cup U}:\ \limsup_{\zeta\to z}
\omega(\zeta,S,D\cup U)\ge\eps \}
\end{equation*}
is a compact connected subset of $\overline{D\cup U}$.
Moreover, if $U_1\subset U_2$ then $F^\eps_{U_1}\subset F^\eps_{U_2}$.
We set $F^\eps=\cap_U F^\eps_{U}$. Then $F^\eps$ is a compact connected
subset of $\overline D$.

Now let $z\in K$ be a regular boundary point of $D$.
Then \begin{equation}
\lim_{\zeta\to z, \zeta\in D}\omega(\zeta,S,D)=0, \label{j1}
\end{equation}
because it is the solution of a Dirichlet problem with continuous boundary
values. Moreover, if $z\notin K$, then for sufficiently small neighborhoods
$U$ of $K$ we find that $z$ is a regular boundary point of $D\cup U$.
We conclude that $F^\eps\cap \partial D$ is a subset of the union of $K$
and the irregular boundary points of $D$, hence $F_\eps\cap \partial D$ is
a polar set and therefore totally disconnected.

To reach a contradiction, suppose that
$$\inf\{\omega(z,S,D\cup U): K\subset U\text{ open}\}=\eps>0,$$
then $z\in F^\eps$. 
Note that for any decreasing sequence $\{U_i\}$ of neighborhoods of $K$
with $K=\cap U_i$, the functions
$\omega(z,S,D\cup U_i)$ form a decreasing sequence of harmonic functions on 
$D\setminus S$ and therefore converge uniformly on compact sets in $D\setminus S$ to $\omega(z,S,D)$. 

In view of (\ref{j1}) we conclude that there is a neighborhood $V$ of $z$
such that $F^\eps\cap V\subset \partial D$. Thus $F^\eps\cap V$ is
totally disconnected and $\{z\}$ is a component of $F^\eps$, which is a
contradiction.
\end{proof}
%%%%%%%%%%%%%%%%%%%%%%%%%%%%%%%%%%%%%%%%%%%%%%%%%%%%%%%%%%%%%%%%%%%%%%

Let $S$ be a subset of $\CC$ and let $\zeta\in\CC$.
Recall (cf.~\cite{R}) that $S$ is called {\em thin} at $\zeta$ if
$\zeta$ is not a point of density of $S$ or there exists  a subharmonic
function $u$ defined in a neighborhood of $\zeta$ with
\begin{equation}
\limsup_{\overset{z\to\zeta}{z\in S\setminus\zeta}}u(z)<u(\zeta).\label{thin}
\end{equation}
The definition of {\em pluri-thin} is entirely analogous: Instead of working
in $\CC$ one works in $\CC^n$, and the subharmonic function $u$ must be
replaced by a plurisubharmonic function.

It is known that it is equivalent to demand in (\ref{thin}) that the limsup
equals $-\infty$, compare the proof of the Wiener criterion in \cite{R},
and for the several variable case \cite{K}.

Next we see that thin boundary points behave as interior points as far as
harmonic measure is concerned. 

\begin{prop}\label{lemma} 
Let $D$ be a domain in $\CC$ and let $B\subset D$ be a closed subset.
Suppose that $z_0\in D$ is such that $B$ is thin at $z_0$.
Then there exist a sequence of radius $\{r_n\}_{n=1}^\infty\subset(0,1)$,
$r_n\to0$, and a sequence $\{z_k\}_{k=1}^\infty\subset D\setminus B$, 
$z_k\to z_0$,  such that  $\DD(z_0,r)\Subset D$, 
$\partial\DD(z_0,r)\subset D\setminus B$, and 
\begin{equation*}
\omega(z_k,\partial\DD(z_0,r_n),\DD(z_0,r_n)\setminus B)\ge
\frac{1}{2},\quad n,k=1,2,\dots
\end{equation*}
\end{prop}

\begin{proof} 
Since $B$ is thin at $z_0$, we can find a subharmonic function $u$ 
defined on a neighborhood of $z_0$ such that 
\begin{equation*}
\limsup_{\overset{z\to z_0}{z\in B\setminus\{z_0\}}}u(z)=-\infty<u(z_0).
\end{equation*} 
Moreover, by scaling and adding a constant,
we can assume that $u(z_0)=-\frac13$ and $u<0$ 
on a small ball $\DD(z_0,\epsilon)$.

It is an immediate consequence of Wieners criterion that there exists a 
sequence of radii $r_n\searrow 0$ such that the circles
$\partial\DD(z_0,r_n)$ do not intersect $B$
(see e.g.~\cite[Theorem 5.4.2.]{R}).
Pick ${n_0}$ so big that $u|_{B\cap\{0<|z-z_0|<r_{n_0}\}}\le -1$.
We assume that $n\ge n_0$.
We use the fact that $D$ is non-thin at $z_0$ to find  a sequence
$\{z_k\}\subset D\setminus B$ tending to $z_0$ such that  $u(z_k)\to u(z_0)$.
We assume that $u(z_k)\ge-\frac12$.
Note that
\begin{equation*}
-\omega(z,\partial\DD(z_0,r_{n}),\DD(z_0,r_{n})\setminus B)+u(z)\le-1,\quad
z\in\DD(z_0,r_{n})\setminus B.
\end{equation*}
In particular,
$\omega(z_k,\partial\DD(z_0,r_{n}),\DD(z_0,r_{n})\setminus B)\ge\frac12$.
\end{proof}

%%%%%%%%%%%%%%%%%%%%%%%%%%%%%%%%%%%%%%%%%%%%%%%%%%%%%%%%%%%%%%%%%%%%%%
\section{Proofs of main results}
%%%%%%%%%%%%%%%%%%%%%%%%%%%%%%%%%%%%%%%%%%%%%%%%%%%%%%%%%%%%%%%%%%%%%%
The next result gives the proof of implication
$(1)\implies(2)$ and the second part of Theorem~\ref{main-thm}.
%%%%%%%%%%%%%%%%%%%%%%%%%%%%%%%%%%%%%%%%%%%%%%%%%%%%%%%%%%%%%%%%%%%%%%
\begin{prop}\label{main-thm3}
Let $D$ be a domain in $\CC$, let $A$ be a closed polar subset of $D$, and
let $f\in\OO(D\setminus A)$. Assume that there exist a $z_0\in A$ and
an $R>0$ such that $\{z\in D\setminus A: |f(z)|\ge R\}$ is thin at $z_0$.
Then there exists a $w_0\in\CC$, $|w_0|\le R$, 
such that $\Gamma_f$ is not pluri-thin at $(z_0,w_0)$ (and, therefore,
$(z_0,w_0)\in (\Gamma_f)^\ast_{D\times\CC}$).
\end{prop}

%%%%%%%%%%%%%%%%%%%%%%%%%%%%%%%%%%%%%%%%%%%%%%%%%%%%%%%%%%%%%%%%%%%%%%

\begin{proof}
We may assume that $z_0=0$. Put 
$B=A\cup\{z\in D\setminus A: |f(z)|\ge R\}$. Note that $B$ is a closed subset
of $D$ thin at $0$. By Proposition~\ref{lemma} there exist a sequence
of radius $\{r_n\}$, $r_n\to0$,
and a sequence $\{z_k\}\subset D\setminus B$ tending to $0$ such that
\begin{equation*}
\omega(z_k,\partial \DD_{r_n}, \DD_{r_n}\setminus B)\ge\frac12.
\end{equation*}
Let $w_0$ be a limit point of $\{f(z_k)\}$; clearly $|w_0|\le R$. 
By passing to a subsequence we can assume that $f(z_k)\to w_0$.

We shall show that $\Gamma_f$ is not pluri-thin at $(0,w_0)$.
Indeed, assume that $\Gamma_f$ is pluri-thin at $(0,w_0)$. Then
there exists a plurisubharmonic function $u\in\psh(\CC^2)$ such that
\begin{equation*}
\limsup_{\Gamma_f\ni(z,w)\to (0,w_0)} u(z,w)=-\infty<u(0,w_0).
\end{equation*}
Consider the subharmonic functions
\begin{equation*}
h_k(z)=u(z-z_k,f(z)+w_0-f(z_k))
\end{equation*}
defined on $D\setminus A$. We have $h_k(z_k)=u(0,w_0)$. 

Put $H=\sup_{|z|\le2, |w|\le 3R} u(z,w)$ and put
$C_{nk}=\sup_{z\in\partial\DD_{r_n}} h_k(z)$, $n,k\in\NN$. Note that
$C_{nk}\to C_{n}=\sup_{z\in\partial\DD_{r_n}} u\big(z,f(z)\big)$,
when $k\to\infty$. Moreover, $C_n\to-\infty$ when $n\to\infty$.

We have $h_k\le H$ for $z\in\DD_{r_n}\setminus B$.
We estimate $u(0,w_0)=h_k(z_k)$ using the two-constants theorem:
\begin{equation*}
u(0,w_0)=h_k(z_k)
\le H-(H-C_{nk})\omega(z_k,\partial\DD_{r_n},\DD_{r_n}\setminus B)
\le \frac{H+C_{nk}}{2}.
\end{equation*}
Now we let $k\to\infty$. So, $u(0,w_0)\le\frac{H+C_n}{2}$. If we let
$n\to\infty$, we get $u(0,w_0)=-\infty$. A contradiction.
\end{proof}
%%%%%%%%%%%%%%%%%%%%%%%%%%%%%%%%%%%%%%%%%%%%%%%%%%%%%%%%%%%%%%%%%%%%%%
For the proof of implication $(2)\implies(1)$ in Theorem~\ref{main-thm}
we need the following technical result, which follows from the localization 
principle (Theorem~\ref{thm-loc}), Theorem~\ref{thm-imp7}, and 
Corollary~\ref{thm10}.
%%%%%%%%%%%%%%%%%%%%%%%%%%%%%%%%%%%%%%%%%%%%%%%%%%%%%%%%%%%%%%%%%%%%%%
\begin{corollary}\label{cor-12}
Let $D$ be a domain in $\CC$ and let $A$ be a closed polar subset of $D$.
Assume that $f\in\OO(D\setminus A)$ and that $S\subset D\setminus A$ is
a closed disc. Then for any domain $G\Subset D$ and any $R_2>R_1>0$ such that
$S\subset G$, $A\cap\partial G=\varnothing$, $R_1>\max_{z\in S} |f(z)|$ we have
\begin{equation*}
\omega(z,S,G_{R_1})\le
\omega\big((z,f(z)),\Gamma_S,G\times\DD_{R_2})
\le\omega(z,S,G_{R_2}),\quad z\in G_{R_1},
\end{equation*}
where $\Gamma_S=\{\big(z,f(z)\big): z\in S\}$ and
$G_R=\{z\in G\setminus A: |f(z)|<R\}$.
Moreover, if $z_0\in A$ is a regular point of $G_{R_2}$ then
\begin{equation*}
\omega\big((z_0,w_0),\Gamma_S,G\times\DD_{R_1})=0,\quad\text{ for any }
w_0\in\DD_{R_1}.
\end{equation*}
\end{corollary}

\begin{proof} 
Note that $Z=\Gamma_f\cup(A\times\CC)$ is a closed complete pluripolar
set in $D\times\CC$ (use Corollary~\ref{thm10}). 

Fix a neighborhood $V$ of $A\cap G$. Put 
$\widetilde V=V\cup(D\setminus \overline{G})$ and put
\begin{multline*}
U=\big(\{z\in D\setminus A: |f(z)|<R_2\}\cup\widetilde V\big)\times\CC\\
\bigcup\{z\in D\setminus A: |f(z)|>R_1\}\times\{w\in\CC:|w|>R_1\}.
\end{multline*}
Note that $U$ is a neighborhood of $Z$. 
According to the localization principle we have
\begin{multline*}
\omega((z,w),\Gamma_S,G\times\DD_{R_1})=
\omega((z,w),\Gamma_S,(G\times\DD_{R_1})\cap U),\\ %\quad
(z,w)\in\big(G\times\DD_{R_1}\big)\cap U.
\end{multline*}
We have $(G\times\DD_{R_1})\cap U=\big(G_{R_2}\cup V\big)\times\DD_{R_1}$, 
hence
\begin{multline*}
\omega\big((z,w),\Gamma_S,G\times\DD_{R_1}\big)=
\omega\big((z,w),\Gamma_S,(G_{R_2}\cup V)\times\DD_{R_1}\big),\\
(z,w)\in\big(G_{R_2}\cup V\big)\times\DD_{R_1}.
\end{multline*}
Now it suffices to use the inequalities
\[
\omega\big((z,f(z)),\Gamma_S,G\times\DD_{R_1}\big)\ge
\omega\big(z,S,G_{R_1}\big),\quad z\in G_{R_1},
\]
and
\begin{multline*}
\omega\big((z,w),\Gamma_S,(G_{R_2}\cup V)\times\DD_{R_1}\big)\le
\omega\big(z,S,G_{R_2}\cup V\big),\\  %\quad 
(z,w)\in \big(G_{R_2}\cup V\big)\times\DD_{R_1},
\end{multline*}
and Theorem~\ref{thm-imp7}.
\end{proof}

%%%%%%%%%%%%%%%%%%%%%%%%%%%%%%%%%%%%%%%%%%%%%%%%%%%%%%%%%%%%%%%%%%%%%%

\begin{proof}[Proof of Theorem~\ref{main-thm}]
Note that all we have to prove is the 
implication $(2)\implies(1)$.

Let $S$ be a closed disc in $D\setminus A$. Put
\begin{equation*}
\Gamma_S:=\{\big(z,f(z)\big):z\in S\}.
\end{equation*}
Let $M$ be the maximum of $|f|$ on $S$ and let $R_\nu\nearrow+\infty$.

Fix domains $D_\nu\Subset D_{\nu+1}\subset D$ such that $\cup_\nu D_\nu=D$
and $\partial D_\nu\cap A=\varnothing$. According to Corollary~\ref{cor-12}
\begin{equation*}
\omega\big((z_0,w_0),\Gamma_S,D_\nu\times\DD_{R_\nu})=0,\quad\text{ for any }
w_0\in\DD_{R_\nu}, \nu\in\NN.
\end{equation*}
Hence, $(\Gamma_S)^{-}_{D_\nu\times R_\nu}\cap (\{z_0\}\times\CC)=\varnothing$,
for any $\nu\in\NN$. So, 
$(\Gamma_S)^\ast_{D\times\CC}\cap (\{z_0\}\times\CC)=\varnothing$.
It suffices to remark that 
$(\Gamma_S)^\ast_{D\times\CC}=(\Gamma_f)^\ast_{D\times\CC}$.
\end{proof}
%%%%%%%%%%%%%%%%%%%%%%%%%%%%%%%%%%%%%%%%%%%%%%%%%%%%%%%%%%%%%%%%%%%%%%

Theorem~\ref{thm11} immediately follows from the following more
general result.
\begin{theorem}[cf.~\cite{W2}]\label{thm1'}
Let $A$ be a closed polar set in a domain $D\subset\CC$ 
and let $f\in\OO(D\setminus A)$. Let $B$ be a subset of $A$
such that $B\cap A'=\varnothing$.
Assume that $f$ does not extend holomorphically to any point of $B$.
Then the set 
\begin{equation*}
\{(z,f(z)):z\in D\setminus A\}\bigcup \big((A\setminus B)\times\CC\big)
\end{equation*}
is complete pluripolar in $D\times\CC$.
\end{theorem}

\begin{proof} 
Note that the set $B$ has only isolated points. 
According to Corollary~\ref{thm10} we know that $\Gamma_f\cup (A\times\CC)$
is complete pluripolar. Since $\{z\in D\setminus A: |f(z)|\ge R\}$
is not-thin at any point of $B$ (because $f$ has isolated singularities
at any point of $B$), by Theorem~\ref{main-thm} we have the proof.
\end{proof}

The  following result gives some information
on the set $(\Gamma_f)^\ast_{D\times\CC}$ in case $f$ is not extandable
over a polar subset $A$.
\begin{theorem} Let $D$ be a domain in $\CC$ and let $A\subset D$ be a closed
polar subset. Assume that $f\in\OO(D\setminus A)$ and $f$ does not
extend holomorphically over $A$. Then there exists a dense subset $B\subset A$
such that $(\Gamma_f)^\ast_{D\times\CC}\cap (B\times\CC)=\varnothing$.
\end{theorem}

\begin{proof} According to a result from \cite{C} there exists a
dense subset $B\subset A$ with the following property:
for any $b\in B$ there exists
a continuous curve $\gamma:(0,1]\to D\setminus A$ such that 
$\gamma(t)\to b$ and 
$|f\big(\gamma(t)\big)|\to\infty$ when $t\to0$. Hence, for any
$b\in B$ and any $R>0$ the set
$\{z\in D\setminus A: |f(z)|\ge R\}$ is not thin at $b$. So, the 
result follows from Theorem~\ref{main-thm}
\end{proof}

\section{Examples}
In case $\Gamma_f\neq(\Gamma_f)^\ast_{D\times\CC}$, 
it is natural to try to determine the non-trivial part of 
$(\Gamma_f)^\ast_{D\times\CC}$. In Proposition~\ref{verylast}
we will compute  the pluripolar hull of particular examples.
We introduce some notation for this section. Let  
$A=\{a_n\}_{n=1}^\infty\subset\DD\setminus\{0\}$ be a sequence 
such that $a_n\to0$ and let $\{c_n\}_{n=1}^\infty\subset\CC\setminus\{0\}$.
Put
\begin{equation*}
f(z)=\sum_{n=1}^\infty\frac{c_n}{z-a_n}.
\end{equation*}
We will suppose that $\sum_{n=1}^\infty|c_n|<+\infty$ and that
$\sum_{n=1}^\infty\left|\frac{c_n}{a_n}\right|$ is convergent.
We also introduce  $\gamma_N=\sum_{n=N}^{\infty}|c_n|$.
Note that $f\in\OO\big(\CC\setminus (A\cup\{0\})\big)$ and that
$f(0)$ is well-defined.
First let us show the following result.
%%%%
%%%%
\begin{prop}\label{lastprop}
Suppose that there exists a sequence $N_m\nearrow\infty$ such that
\begin{equation}\label{eq-ex}
\lim_{m\to\infty}
\frac{\sum_{n=1}^{N_m}\log|a_n|}
{\log\gamma_{N_m+1}}=0.
\end{equation}
Then 
\[
(\Gamma_f)^\ast_{\CC^2}\subset\Gamma_f\cup \big\{(0,f(0))\big\}.
\]
\end{prop}

\begin{proof} Since $a_n\to0$, 
$\lim_{N\to\infty}\frac{1}{N}\sum_{n=1}^N\log|a_n|=-\infty$.
Hence, from \eqref{eq-ex} we have
$\lim_{m\to\infty}\frac{1}{N_m}\log\gamma_{N_m+1}=-\infty$.
Put 
\[
f_N(z)=\sum_{n=1}^{N}\frac{c_n}{z-a_n}
\]
and
\[
h_N(z,w)=\frac{1}{N}\log\Big|\big(w-f_N(z)\big)\prod_{n=1}^N(z-a_n)\Big|.
\]
Note that $h_N$ is a plurisubharmonic function on $\CC^2$.
Fix $R>0$ and let $|z|,|w|<R$. We have
\begin{multline*}
h_N(z,w)\le\frac{1}{N}\log\Big(
R\prod_{n=1}^N(R+|a_n|)+
\sum_{k=1}^{N}\frac{|c_k|}{R+|a_k|}\prod_{n=1}^N(R+|a_n|)\Big)\\
\le
\frac{1}{N}\log\Big(R+\sum_{k=1}^{N}\frac{|c_k|}{R}\Big)+
\frac{1}{N}\log\Big(\prod_{n=1}^N(R+|a_n|)\Big):=M_N.
\end{multline*}
Observe that $M_N\to\log R$ when $N\to\infty$.

Fix a closed disc $\overline{\DD}(z_0,2r_0)\subset\CC\setminus\overline{A}$.
Let us estimate $h_N\big(z,f(z)\big)$ for $z\in\DD(z_0,r_0)$.
We have
\begin{equation*}
h_N\big(z,f(z)\big)\le\frac{1}{N}\log\Big(
\left(\sum_{n=N+1}^{\infty}\frac{|c_n|}{r_0}\right)\prod_{n=1}^N(r_0+|a_n|)\Big):=K_N.
\end{equation*}
Note that 
\[
\lim_{m\to\infty}
\frac{\frac{1}{N_m}\log\gamma_{N_m+1}}{K_{N_m}}=1.
\]
Put $v_N(z,w)=\frac{M_N-h_N(z,w)}{M_N-K_N}$. If $w\ne f(0)$ then
\[
\lim_{m\to\infty} v_{N_m}(0,w)=
\lim_{m\to\infty}
\frac{\sum_{n=1}^{N_m}\log|a_n|}{\log\gamma_{N_m+1}}=0.
\]
From the definition we have
\begin{equation*}
\omega((z,w),\Gamma_S,\DD_R\times\DD_R)\le v_N(z,w),\quad z,w\in\DD_R,
\end{equation*}
where $S=\overline{\DD}(z_0,r_0)$ and $\Gamma_S=\{(z,f(z)): z\in S\}$.
So, if $w\ne f(0)$ and $w\in\DD_R$ then
\begin{equation*}
0\le\omega((0,w),\Gamma_S,\DD_R\times\DD_R)\le \lim_{m\to\infty}
v_{N_m}(0,w)=0.
\end{equation*}
Hence, $(0,w)\not\in (\Gamma_f)^\ast_{\CC^2}$.
\end{proof}

\begin{prop}\label{verylast}
Suppose that $\sum_{n=1}^\infty\frac{\log|a_n|}{\log \gamma_n}<+\infty$.
Then 
\[
(\Gamma_f)^\ast_{\CC^2}=\Gamma_f\cup \big\{(0,f(0))\big\}.
\]
\end{prop}

\begin{proof}
Fix $N>k>1$. Note that 
\[
\frac{\sum_{n=1}^N\log|a_n|}{\log\gamma_{N+1}}\le
\frac{\sum_{n=1}^k\log|a_n|}{\log\gamma_{N+1}}+
\sum_{n=k+1}^N\frac{\log|a_n|}{\log\gamma_n}.
\]
From this inequality we get that condition \eqref{eq-ex} is satisfied and 
by Proposition~\ref{lastprop} we have 
$(\Gamma_f)^\ast_{\CC^2}\subset\Gamma_f\cup\big\{\big(0,f(0)\big)\big\}$.

So, it suffices to show that $(\Gamma_f)^\ast_{\CC^2}\ne \Gamma_f$.
It is not difficult to see that $\frac{|a_n|}{\sqrt \gamma_n}\to\infty$,
when $n\to\infty$. Hence, there exists a constant $C>0$ such that
$|a_n|>r_n$ for any $n\in\NN$, where $r_n=C\sqrt{\gamma_n}$. Because $\sum_{n=1}^\infty \frac1{|\log\gamma_n|}<\infty$,
we find
\[
R:=\sum_{n=1}^\infty\frac{|c_n|}{r_n} \le \sum_{n=1}^\infty\sqrt{\gamma_n}\frac{|c_n|}{C\gamma_n}=
\frac1C\sum_{n=1}^\infty\sqrt{\gamma_n}<+\infty.
\]
Put $A=\{0\}\cup\{a_n: n\in\NN\}$ and put
$F:=\cup_{n=1}^\infty\overline{\DD(a_n,r_n)}$.
We have
\[
\{z\in\CC\setminus A:|f(z)|\ge R\}\subset F.
\]
According to Theorem~\ref{main-thm}, it suffices to show that $F$ is 
thin at $0$. But this is a variation on Wiener's
criterion, see e.g.~\cite{R}. Indeed, 
there exists a sequence $\alpha_n\nearrow\infty$ such that
\[
\sum_{n=1}^\infty \alpha_n\frac{\log\Big(\frac{|a_n|}{1+|a_n|}\Big)}{\log r_n}
<+\infty.
\]
Consider the function
\[
u(z)=\sum_{n=1}^\infty\frac{\alpha_n}{\log\frac1{r_n}}
\big[\log|z-a_n|-\log(1+|a_n|)\big].
\]
Note that $u$ is a negative subharmonic function on $\DD$ and
that $u(0)>-\infty$. We have 
\[
u(z)\le\frac{\alpha_n}{\log\frac1{r_n}}\big[\log r_n-\log(1+|a_n|)\big],
\quad z\in\overline{\DD(a_n,r_n)}.
\]
Hence,
\[
\limsup_{F\ni z\to 0}
u(z)\le\limsup_{n\to\infty}
\frac{\alpha_n}{\log\frac1{r_n}}\big[\log r_n-\log(1+|a_n|)\big]=
-\infty<u(0).
\]
We obtain that $F$ is thin at $0$.
\end{proof}

\begin{example}\label{example} Take $a_n=\frac{1}{n}$ and $c_n=e^{-n^2}/n^2$,
$n\in\NN$. Put
\[
f(z)=\sum_{n=1}^\infty \frac{c_n}{z-a_n}.
\]
Then by Propostion~\ref{verylast} we have
$(\Gamma_f)^\ast_{\CC^2}=\Gamma_f\cup\big\{\big(0,f(0)\big)\big\}$.
\end{example}

%%%%%%%%%%%%%%%%%%%%%%%%%%%%%%%%%%%%%%%%%%%%%%%%%%%%%%%%%%%%%%%%%%%%%%
%
%%%%%%%%%%%%%%%%%%%%%%%%%%%%%%%%%%%%%%%%%%%%%%%%%%%%%%%%%%%%%%%%%%%%%%
\bibliographystyle{amsplain}

%%%%%%%%%%%%%%%%%%%%%%%%%%%%%%%%%%%%%%%%%%%%%%%%%%%%%%%%%%%%%%%%%%%%%%

\end{document}